\documentclass[12pt]{article}
\usepackage{amsmath,amssymb,amsthm}

\newtheorem{lemma}{Lemma}[section]

\newtheorem{proposition}{Proposition}[section]
\newtheorem{theorem}{Theorem}[section]
\newtheorem{corollary}{Corollary}[section]
\newtheorem{remark}{Remark}[section]

\theoremstyle{definition}

\begin{document}
\title{\LARGE \bf On Besicovitch almost periodic selections of multivalued maps}
\author{\Large L.I.~Danilov \medskip \\
\large Physical-Technical Institute \\
\large Russia, 426000, Izhevsk, Kirov st., 132 \\
\large e-mail: danilov@otf.pti.udm.ru}
\date{}
\maketitle

\begin{abstract}
We prove that Besicovitch almost periodic multivalued maps ${\bf R}\ni t\to
F(t)\in {\rm cl}\, {\mathcal U}$ have Besicovitch almost periodic selections, where
${\rm cl}\, {\mathcal U}$ is the collection of non-empty closed sets of a 
complete metric space ${\mathcal U}$. \par
\end{abstract}

\large

\section*{Introduction}
The objective of the present paper is to prove the existence of Besicovitch
almost periodic (a.p.) selections of Besicovitch a.p. multivalued maps. The
existence of Stepanov a.p. selections of Stepanov a.p. multivalued maps was
proved in \cite{1}. The Stepanov a.p. selections which satisfy some additional
conditions were studied in \cite{2,3,4,5}. The papers \cite{6,7} were devoted
to investigation of Weyl a.p. selections of Weyl a.p. multivalued maps.

In proofs suggested in the paper we use technique from \cite{2,3,6}.

The results of this paper are applied to the study of a.p. solutions of
differential inclusions \cite{8,9}.

In Section 1 we present some properties of Besicovitch a.p. functions which
will be used in what follows (as regards definitions and assertions on a.p.
functions, see e.g. \cite{10}). The main results are contained in Section 2.
We prove the Theorem \ref{th2.1} from Section 2 in Section 3, and the Theorem 
\ref{th3.1} from Section 3 in Section 4.

\section{Some propeties of Besicovitch a.p. functions}

Let $({\mathcal U},\rho )$ be a complete metric space, $\overline A$ the closure
of a set $A\subseteq {\mathcal U}$, $U_r(x)=\{ y\in {\mathcal U}:\rho (x,y)<r\} $, 
$x\in {\mathcal U}$, $r>0$. Let ${\rm meas}$ be Lebesgue measure on ${\bf R}$.
A function $f:{\bf R}\to {\mathcal U}$ is said to be {\it elementary} if there 
exist points $x_j\in {\mathcal U}$ and disjoint measurable (in the Lebesgue sense) 
sets $T_j\subseteq {\bf R}$, $j\in {\bf N}$, such that ${\rm meas}\, {\bf R}\backslash 
\bigcup\limits_jT_j=0$ and $f(t)=x_j$ for all $t\in T_j\, $. We denote this function 
by $f(.)=\sum\limits_jx_j\chi _{T_j}(.)$ (where $\chi _T(.)$ is the characteristic 
function of a set $T\subseteq {\bf R}$). For arbitrary functions $f_j:{\bf R}\to 
{\mathcal U}$, $j\in {\bf N}$, we define the function $\sum\limits_jf_j(.)\chi 
_{T_j}(.):{\bf R}\to {\mathcal U}$ that coincides with functions $f_j(.)$ on the sets 
$T_j\, $, $j\in {\bf N}$ (the notation $\sum\limits_jf_j(.)\chi _{T_j}(.)$ will be
used not only in the case when ${\mathcal U}=({\mathcal H},\| .\| )$ is Banach
space but also for arbitrary metric space ${\mathcal U}=({\mathcal U},\rho )$, and
in last case no linear operations will be carried out on the functions under 
consideration). A function $f:{\bf R}\to {\mathcal U}$ is {\it measurable}
if for any $\epsilon >0$ there exists an elementary function $f_{\epsilon}:{\bf 
R}\to {\mathcal U}$ such that
$$
{\rm ess}\, \sup\limits_{\hskip -0.7cm t\in {\bf R}}\, \rho (f(t),f_{\epsilon}(t))
<\epsilon \, .
$$

Let $M({\bf R},{\mathcal U})$ be the set of measurable functions $f:{\bf R}\to
{\mathcal U}$ (functions that coincide for a.e. $t\in {\bf R}$ will be identified),
$(L^{\infty}({\bf R},{\mathcal U}),D_{\infty})$ the space of essentially bounded
functions from $M({\bf R},{\mathcal U})$ with metric
$$
D_{\infty}(f,g)={\rm ess}\, \sup\limits_{\hskip -0.7cm t\in {\bf R}}\, \rho 
(f(t),g(t))\, ,\ f,g\in L^{\infty}({\bf R},{\mathcal U})\, .
$$
Let a point $x_0\in {\mathcal U}$ be fixed. We use the notation
$$
M_p({\bf R},{\mathcal U})\doteq \{ f\in M({\bf R},{\mathcal U}):\sup\limits_{\xi
\in {\bf R}}\, \int\limits_{\xi}^{\xi +1}{\rho }^{\, p}(f(t),x_0)\, dt< +\infty 
\} \, ,\ p\geq 1\, ,
$$
and define the metric on $M_p({\bf R},{\mathcal U})$:
$$
D^{(S)}_p(f,g)=\biggl( \, \sup\limits_{\xi \in {\bf R}}\, 
\int\limits_{\xi}^{\xi +1}{\rho }^{\, p}(f(t),g(t))\, dt\biggr) ^{1/p},\ f,g\in
M_p({\bf R},{\mathcal U})\, .
$$

For a Banach space ${\mathcal U}=({\mathcal H},\| .\| )$ ($\rho (x,y)=
\| x-y\| $, $x,y\in {\mathcal H}$; $\| x\| =|x|$ if $x\in {\bf R}$) we 
denote by
$$
\| f\| _{\infty}={\rm ess}\, \sup\limits_{\hskip -0.7cm t\in {\bf R}}\,
\| f(t)\| \, ,\ f\in L^{\infty}({\bf R},{\mathcal H})\, ,
$$
and 
$$
\| f\| ^{(S)}_p=\biggl( \, \sup\limits_{\xi \in {\bf R}}\, 
\int\limits_{\xi}^{\xi +1}\| f(t)\| ^pdt\biggr) ^{1/p},\ f\in
M_p({\bf R},{\mathcal H})\, ,
$$
the norms on linear spaces $L^{\infty}({\bf R},{\mathcal H})$ and $M_p({\bf 
R},{\mathcal H})$, $p\geq 1$, respectively. In what follows, we shall use
the notation ${\mathcal H}$ for Banach space, and it will be convenient 
to assume the Banach space ${\mathcal H}=({\mathcal H},\| .\| )$ to be complex. 
If the Banach space ${\mathcal H}$ is real, then we can consider the 
complexification ${\mathcal H}+i{\mathcal H}$ identifying the space ${\mathcal H}$ 
with the real subspace (the norm $\| .\| _{{\mathcal H}+i{\mathcal H}}$ on the 
real subspace coincides with the norm $\| .\| $).

A set $T\subseteq {\bf R}$ is called {\it relatively dense} if there exists a
number $a>0$ such that $[\xi ,\xi +a]\cap T\neq \emptyset $ for all $\xi \in 
{\bf R}$. A number $\tau \in {\bf R}$ is called an {\it 
$(\epsilon ,D_{\infty})$-almost period} of a function $f\in L^{\infty}({\bf 
R},{\mathcal U})$, $\epsilon >0$, if $D_{\infty}(f(.),f(.+\tau ))<\epsilon $.
A continuous function $f\in C({\bf R},{\mathcal U})\cap L^{\infty}({\bf 
R},{\mathcal U})$ belongs to the space $CAP\, ({\bf R},{\mathcal U})$ of 
{\it Bohr a.p.} functions if for any $\epsilon >0$ the set of $(\epsilon ,
D_{\infty})$-almost periods of the function $f$ is relatively dense. A number
$\tau \in {\bf R}$ is called an {\it $(\epsilon ,D^{(S)}_p)$-almost period} 
of a function $f\in M_p({\bf R},{\mathcal U})$, $p\geq 1$, if $D^{(S)}_p(f(.),
f(.+\tau ))<\epsilon $. A function $f\in M_p({\bf R},{\mathcal U})$, $p\geq 1$, 
belongs to the space $S_p({\bf R},{\mathcal U})$ of {\it Stepanov a.p.} functions 
{\it of order $p\geq 1$} if for any $\epsilon >0$ the set of $(\epsilon ,
D^{(S)}_p)$-almost periods of $f$ is relatively dense. 

On the space ${\mathcal U}$ we also consider the metric $\rho ^{\, \prime}(x,y)=
\min \, \{ 1,\rho (x,y)\} $, $x,y\in {\mathcal U}$. The metric space $({\mathcal 
U},\rho ^{\, \prime})$ is complete (as well as $({\mathcal U},\rho )$). We
define the metric on $M({\bf R},{\mathcal U})=M_1({\bf R},({\mathcal U},\rho ^{\, 
\prime}))\, $: 
$$
D^{(S)}(f,g)=\sup\limits_{\xi \in {\bf R}}\, \int\limits_{\xi}^{\xi +1}{\rho }^{\, 
\prime}(f(t),g(t))\, dt\, ,\ f,g\in M({\bf R},{\mathcal U})\, .
$$ 
Let $S({\bf R},{\mathcal U})\doteq S_1({\bf R},({\mathcal U},\rho ^{\, 
\prime}))$ ({\it Stepanov a.p.} function $f\in S({\bf R},{\mathcal U})$ is 
defined as Stepanov a.p. function of order 1 taking values in the metric
space $({\mathcal U},\rho ^{\, \prime})$). We have $CAP\, ({\bf R},{\mathcal U})
\subseteq S_p({\bf R},{\mathcal U})\subseteq S_1({\bf R},{\mathcal U})\subseteq 
S({\bf R},{\mathcal U})$.

A sequence $\tau _j\in {\bf R}$, $j\in {\bf N}$, is said to be {\it 
$f$-returning} for a function $f\in S({\bf R},{\mathcal U})$ if 
$D^{(S)}(f(.),f(.+\tau _j))\to 0$ as $j\to +\infty $. If $f
\in CAP\, ({\bf R},{\mathcal U})\subseteq S({\bf R},{\mathcal U})$, then a
sequence $\tau _j\in {\bf R}$, $j\in {\bf N}$, is $f$-returning if and only if
$D_{\infty}(f(.),f(.+\tau _j))\to 0$ as $j\to +\infty $. If $f\in S_p({\bf R},
{\mathcal U})\subseteq S({\bf R},{\mathcal U})$, $p\geq 1$, then a sequence 
$\tau _j\in {\bf R}$, $j\in {\bf N}$, is $f$-returning if and only if 
$D^{(S)}_p(f(.),f(.+\tau _j))\to 0$ as $j\to +\infty $. 

For a function $f\in S({\bf R},{\mathcal U})$ we denote by ${\rm Mod}\, f$ 
the set of numbers $\lambda \in {\bf R}$ for which $e^{\, i\lambda \tau _j}\to 1$
($i^2=-1$) as $j\to +\infty $ for all $f$-returning sequences $\tau _j\, $. The
set ${\rm Mod}\, f$ is a module (additive group) in ${\bf R}$. If a function
$f\in S({\bf R},{\mathcal U})$ is not a.e. (almost everywhere) constant, then
${\rm Mod}\, f$ is a countable module (${\rm Mod}\, f=\{ 0\} $ otherwise). If
${\mathcal U}=({\mathcal H},\| .\| )$ is a Banach space, then for all functions
$f\in S_1({\bf R},{\mathcal H})$ the sets ${\rm Mod}\, f$ coincide with the
modules of Fourier exponents of the functions $f$.

For any function $f\in S_p({\bf R},{\mathcal H})$ and any $\epsilon >0$ there
exists a function $f_{\epsilon}\in CAP\, ({\bf R},{\mathcal H})$ such that 
$\| f-f_{\epsilon}\| ^{(S)}_p<\epsilon $ and ${\rm Mod}\, f_{\epsilon}\subseteq
{\rm Mod}\, f$ (moreover, the Fourier exponents of a function $f_{\epsilon}$
belong to the set of Fourier exponents of a function $f$). If $f\in S({\bf R},
{\mathcal H})$, then for any $\epsilon >0$ there is a function $f_{\epsilon}\in 
CAP\, ({\bf R},{\mathcal H})$ such that $D^{(S)}(f,f_{\epsilon})<\epsilon $ and 
${\rm Mod}\, f_{\epsilon}\subseteq {\rm Mod}\, f$.

Let ${\mathcal M}_p({\bf R},{\mathcal U})$, $p\geq 1$, be the Marcinkiewicz
space, i.e. the set of functions $f\in M({\bf R},{\mathcal U})$ for which
$\rho (f(.),x_0)\in L^p_{\rm loc}({\bf R},{\bf R})$ and
$$
\overline {\lim\limits_{b\to +\infty}}\ \, \frac 1{2b}\ \int\limits_{-b}^b
\rho ^{\, p}(f(t),x_0)\, dt<+\infty \, .
$$
We set
$$
D^{(B)}_p(f,g)=\biggl( \ \overline {\lim\limits_{b\to +\infty}}\ \, \frac 1{2b}\ 
\int\limits_{-b}^b\rho ^{\, p}(f(t),g(t))\, dt\, \biggr) ^{1/p},\ f,g\in {\mathcal 
M}_p({\bf R},{\mathcal U})\, .
$$
If ${\mathcal U}=({\mathcal H},\| .\| )$ is a Banach space, then we define the
seminorm
$$
\| f\| ^{(B)}_p=\biggl( \ \overline {\lim\limits_{b\to +\infty}}\ \, \frac 1{2b}\ 
\int\limits_{-b}^b\| f(t)\| ^p\, dt\, \biggr) ^{1/p},\ f\in {\mathcal M}_p({\bf R},
{\mathcal H})\, .
$$
For functions $f,g\in {\mathcal M}_p({\bf R},{\mathcal U})$ let us define the
equivalence relation: $f\sim g$ if and only if $D^{(B)}_p(f,g)=0$. Then the
quotient space $({\mathcal M}_p({\bf R},{\mathcal U})/\sim ,D^{(B)}_p)$ is
complete metric space \cite{11}. We have $M_p({\bf R},{\mathcal U})
\subseteq {\mathcal M}_p({\bf R},{\mathcal U})$ and $D^{(B)}_p(f,g)\leq
D^{(S)}_p(f,g)$ for all functions $f,g\in M_p({\bf R},{\mathcal U})$.

A function $f\in {\mathcal M}_p({\bf R},{\mathcal U})$, $p\geq 1$, 
belongs to the space $B_p({\bf R},{\mathcal U})$ of {\it Besicovitch a.p.} 
functions {\it of order $p$} if for any $\epsilon >0$ there is a function 
$f_{\epsilon}\in S_p({\bf R},{\mathcal U})$ such that $D^{(B)}_p(f,f_{\epsilon})
<\epsilon $.

By the Fr\'echet Theorem \cite{12} the metric space $({\mathcal U},\rho )$ can
be isometrically embedded into some Banach space ${\mathcal H}$, hence the
following definition of the space $B_p({\bf R},{\mathcal U})$ is equivalent
to the previous one: a function $f\in {\mathcal M}_p({\bf R},{\mathcal U})$
belongs to the space $B_p({\bf R},{\mathcal U})$, if for some Banach space
${\mathcal H}$ into which the metric space $({\mathcal U},\rho )$ is
isometrically embedded (and therefore for all such Banach spaces ${\mathcal H}$)
and for all $\epsilon >0$ there exists a function $f_{\epsilon}\in CAP\, ({\bf 
R},{\mathcal H})$ such that $\|f-f_{\epsilon}\| ^{(B)}_p<\epsilon $ (where $\| .
\| ^{(B)}_p$ is the seminorm on the space ${\mathcal M}_p({\bf R},{\mathcal H})$
and we assume that the function $f$ takes values in the space ${\mathcal H}$).

For functions $f,g\in M({\bf R},{\mathcal U})={\mathcal M}_1({\bf R},({\mathcal U},
\rho ^{\, \prime}))$ we denote
$$
D^{(B)}(f,g)=\, \overline {\lim\limits_{b\to +\infty}}\ \, \frac 1{2b}\ 
\int\limits_{-b}^b\rho ^{\, \prime}(f(t),g(t))\, dt\, .
$$
Let $B({\bf R},{\mathcal U})\doteq B_1({\bf R},({\mathcal U},\rho ^{\, 
\prime}))$ be the space of {\it Besicovitch a.p.} functions (defined as
Besicovitch a.p. functions of order 1 taking values in the metric space 
$({\mathcal U},\rho ^{\, \prime})$). We have $S({\bf R},{\mathcal U})
\subseteq B({\bf R},{\mathcal U})$ and $S_p({\bf R},{\mathcal U})\subseteq 
B_p({\bf R},{\mathcal U})\subseteq B_1({\bf R},{\mathcal U})\subseteq B({\bf 
R},{\mathcal U})$.

A sequence $\tau _j\in {\bf R}$, $j\in {\bf N}$, is said to be {\it 
$f$-returning} for a function $f\in B({\bf R},{\mathcal U})$ if 
$D^{(B)}(f(.),f(.+\tau _j))\to 0$ as $j\to +\infty $. If $f
\in S({\bf R},{\mathcal U})\subseteq B({\bf R},{\mathcal U})$, then a
sequence $\tau _j\in {\bf R}$, $j\in {\bf N}$, is $f$-returning if and only if
$D^{(S)}(f(.),f(.+\tau _j))\to 0$ as $j\to +\infty $. If $f\in B_p({\bf R},
{\mathcal U})\subseteq B({\bf R},{\mathcal U})$, $p\geq 1$, then a sequence 
$\tau _j\in {\bf R}$, $j\in {\bf N}$, is $f$-returning if and only if 
$D^{(B)}_p(f(.),f(.+\tau _j))\to 0$ as $j\to +\infty $. (The set of 
$f$-returning sequences is determined only by the a.p. function itself 
and does not depend on the spaces under consideration of a.p. functions
which include the function $f$.)

For a function $f\in B({\bf R},{\mathcal U})$ (by analogy with a function
$f\in S({\bf R},{\mathcal U})$) we denote by ${\rm Mod}\, f$ the set 
(module) of numbers $\lambda \in {\bf R}$ for which $e^{\, i\lambda \tau _j}\to 1$
as $j\to +\infty $ for all $f$-returning sequences $\tau _j\, $. If there
exists a constant function $y(t)\equiv y\in {\mathcal U}$, $t\in {\bf R}$,
such that $D^{(B)}(f(.),y(.))=0$, then ${\rm Mod}\, f=\{ 0\} $. If
$D^{(B)}(f(.),y(.))\neq 0$ for all constant functions $y(t)\equiv y\in {\mathcal 
U}$, $t\in {\bf R}$, then ${\rm Mod}\, f$ is a countable module.

If $f\in B({\bf R},{\mathcal U})$ and $\tau _j\in {\bf R}$, $j\in {\bf N}$,
is a sequence for which $e^{\, i\lambda \tau _j}\to 1$ as $j\to +\infty $ for all
numbers $\lambda \in {\rm Mod}\, f$, then $\tau _j$ is $f$-returning sequence.

For a function $f\in B_1({\bf R},{\mathcal H})$ we denote by $\Lambda \{ f\} $
the set of its Fourier exponents, i.e. the set of numbers $\lambda \in {\bf R}$
for which
$$
\lim\limits_{b\to +\infty}\ \, \frac 1{2b}\ \int\limits_{-b}^be^{-i\lambda t}
f(t)\, dt\neq 0
$$
(the limit exists for all numbers $\lambda \in {\bf R}$). The module ${\rm Mod}\, 
f$ of a function $f\in B_1({\bf R},{\mathcal H})$ coincides with the module
of Fourier exponents $\lambda \in \Lambda \{ f\} $, i.e. the smallest module
(additive group) in ${\bf R}$ including the set $\Lambda \{ f\} $.

If $\Lambda _j\subseteq {\bf R}$ are arbitrary modules (the set of indices
$j$ may be an arbitrary non-empty set), then by $\sum\limits_j \Lambda _j $ 
(or by $\Lambda _1+\dots +\Lambda _n$ for finitely many modules 
$\Lambda _j\, $, $j=1,\dots ,n$) we denote the sum of modules, that is, the
smallest module (additive group) in ${\bf R}$ containing all the sets $\Lambda 
_j\, $.

Suppose that $f\in B({\bf R},{\mathcal U})$ and $f_j\in B({\bf R},{\mathcal 
U}_j)$, $j\in {\bf N}$, where the ${\mathcal U}_j$ are (complete) metric spaces.
Then ${\rm Mod}\, f\subseteq \sum\limits_j{\rm Mod}\, f_j$ if and only if every 
sequence $\tau _k\in {\bf R}$, $k\in {\bf N}$, which is $f_j$-returning for all 
$j\in {\bf N}$ is $f$-returning. In particular, if $f_j\in B({\bf R},
{\mathcal U}_j)$, $j=1,2$, then ${\rm Mod}\, f_1\subseteq {\rm Mod}\, f_2$ if 
and only if every $f_2$-returning sequence $\tau _k\in {\bf R}$, $k\in {\bf N}$, is 
$f_1$-returning.

If $f\in M({\bf R},{\mathcal U})$, $f_j\in B({\bf R},{\mathcal U})$, $j\in 
{\bf N}$, and $D^{(B)}(f,f_j)\to 0$ as $j\to +\infty $, then $f\in B({\bf R},
{\mathcal U})$ and ${\rm Mod}\, f\subseteq \sum\limits_j{\rm Mod}\, f_j\, $.

\begin{proposition} \label{p1.1}
For every function $f\in B_p({\bf R},{\mathcal H})$, $p\geq 1$ (where ${\mathcal
H}$ is a complex Banach space), and every $\epsilon >0$ there
exists a function $f_{\epsilon}\in CAP\, ({\bf R},{\mathcal H})$ such that 
$\| f-f_{\epsilon}\| ^{(B)}_p<\epsilon $ and $\Lambda \{ f_{\epsilon}\} \subseteq
\Lambda \{ f\} $. If $f\in B({\bf R},{\mathcal H})$, then for every $\epsilon >0$ 
there exists a function $f_{\epsilon}\in CAP\, ({\bf R},{\mathcal H})$ such that 
$D^{(B)}(f,f_{\epsilon})<\epsilon $ and $\Lambda \{ f_{\epsilon}\} \subseteq 
{\rm Mod}\, f$.
\end{proposition}

\begin{proposition} \label{p1.2}
For every function $f\in B_p({\bf R},{\mathcal U})$, $p\geq 1$, and every 
$\epsilon >0$ there exists a function $f_{\epsilon}\in S_1({\bf R},{\mathcal 
U})\cap L^{\infty}({\bf R},{\mathcal U})\subseteq S_p({\bf R},{\mathcal U})$ 
such that $D^{(B)}_p(f,f_{\epsilon})<\epsilon $ and ${\rm Mod}\, f_{\epsilon}
\subseteq {\rm Mod}\,  f$. If $f\in B({\bf R},{\mathcal U})$, then for every 
$\epsilon >0$ there exists a function $f_{\epsilon}\in S_1({\bf R},{\mathcal U})
\cap L^{\infty}({\bf R},{\mathcal U})$ such that $D^{(B)}(f,f_{\epsilon})<
\epsilon $ and ${\rm Mod}\, f_{\epsilon} \subseteq {\rm Mod}\, f$.
\end{proposition}

\begin{lemma} \label{l1.1}
Let $({\mathcal U},\rho )$ and $({\mathcal V},\rho _{\mathcal V})$ be (complete) 
metric spaces. Suppose that for a function ${\mathcal F}:{\mathcal U}\to
{\mathcal V}$ there exists a number $C\geq 0$ such that the inequality
$$
\rho _{\mathcal V}({\mathcal F}(u_1),{\mathcal F}(u_2))\leq C\, \rho (u_1,u_2)
$$
holds for all $u_1,u_2\in {\mathcal U}$. Then for all functions $f\in 
B_p({\bf R},{\mathcal U})$ we have ${\mathcal F}(f(.))\in B_p({\bf R},{\mathcal 
V})$ and ${\rm Mod}\, {\mathcal F}(f(.))\subseteq {\rm Mod}\, f(.)$. If $f\in 
B({\bf R},{\mathcal U})$, then also ${\mathcal F}(f(.))\in B({\bf R},{\mathcal 
V})$ and ${\rm Mod}\, {\mathcal F}(f(.))\subseteq {\rm Mod}\, f(.)$.
\end{lemma}

\begin{corollary} \label{c1.1}
Let $f\in B({\bf R},{\mathcal U})$ and $x\in {\mathcal U}$. Then $\rho (f(.),x)
\in B({\bf R},{\bf R})$ and ${\rm Mod}\, \rho (f(.),x)\subseteq {\rm Mod}\, f(.)$.
\end{corollary}

For Banach spaces $({\mathcal H},\| .\| )$ and numbers $a>0$ we define the
functions
$$
{\mathcal H}\ni h\to {\mathcal F}^{\, a}_{\mathcal H}(h)=\left\{
\begin{array}{ll}
h\, ,  &{\text {if}}\ \| h\| \leq a\, ,\\ [0.2cm]
\frac {ah}{\| h\| }\, ,   &{\text {if}}\ \| h\| >a\, .
\end{array}
\right.
$$
For all $h_1,h_2\in {\mathcal H}$ we have
$$
\| {\mathcal F}^{\, a}_{\mathcal H}(h_1)-{\mathcal F}^{\, a}_{\mathcal H}(h_2)\| \leq
2\, \| h_1-h_2\| \, ,
$$
hence the following Lemma \ref{l1.2} is a consequence of Lemma \ref{l1.1}.

\begin{lemma} \label{l1.2}
Let $f\in B({\bf R},{\mathcal H})$. Then for any $a>0$ the function ${\mathcal 
F}^{\, a}_{\mathcal H}(f(.))$ belongs to the set $B({\bf R},{\mathcal H})\cap 
L^{\infty}({\bf R},{\mathcal H})\subset B_1({\bf R},{\mathcal H})$ and ${\rm
Mod}\, {\mathcal F}^{\, a}_{\mathcal H}(f(.))\subseteq {\rm Mod}\, f(.)$.
\end{lemma}

For measurable set $T\subseteq {\bf R}$ let us denote
$$
\widetilde {\kappa}\, (T)=\, \overline {\lim\limits_{b\to +\infty}}\ \, \frac 1{2b}\
{\rm meas}\, [-b,b]\backslash T\, .
$$
For all measurable sets $T_1,T_2\subseteq {\bf R}$ we have $\widetilde {\kappa}\,
(T_1\cap T_2)\leq \widetilde {\kappa}\, (T_1)+\widetilde {\kappa}\, (T_2)$.

Let $f,g\in M({\bf R},{\mathcal U})$, $\epsilon \in (0,1]$ and $\delta >0$.
If $\widetilde {\kappa}\, (\{ t\in {\bf R}:\rho (f(t),g(t))\leq \epsilon \} )
<\delta $, then $D^{(B)}(f,g)\leq \epsilon +\delta $. If $D^{(B)}(f,g)\leq 
\epsilon \delta $, then $\widetilde {\kappa}\, (\{ t\in {\bf R}:\rho (f(t),g(t))
\leq \epsilon \} )\leq {\epsilon}^{-1}D^{(B)}(f,g)\leq \delta $. Hence (see also
Proposition \ref{p1.1}) the Lemma \ref{l1.3} is valid.

\begin{lemma} \label{l1.3}
For any function $f\in B({\bf R},{\mathcal H})$ and any numbers $\epsilon ,
\delta >0$ there is a function $f_{\epsilon ,\, \delta}\in CAP\, ({\bf R},
{\mathcal H})$ such that 
$$
\widetilde {\kappa}\, (\{ t\in {\bf R}:\| f(t)-f_{\epsilon ,\, \delta}(t)\| <
\epsilon \} )<\delta
$$
and $\Lambda \{ f_{\epsilon ,\, \delta}\} \subseteq {\rm Mod}\, f$.
\end{lemma}

The following Lemma \ref{l1.4} is a consequence of Lemma \ref{l1.3} and the
Fr\'echet Theorem.

\begin{lemma} \label{l1.4}
Let $f\in B({\bf R},{\mathcal U})$. Then $\widetilde {\kappa}\, (\{ t\in {\bf R}:
\rho (f(t),x_0)\leq a\} )\to 0$ as $a\to +\infty $.
\end{lemma}

\begin{lemma} \label{l1.5}
Let $f\in B_p({\bf R},{\mathcal U})$, $p\geq 1$. Then
$$
\overline {\lim\limits_{b\to +\infty}}\, \ \frac 1{2b}\ \int\limits_{\{ \, t\, \in
\, [-b,b]\, :\, \rho (f(t),x_0)\, >\, a\, \} }\ \rho ^{\, p}(f(t),x_0)\, dt \to 0
$$
as $a\to +\infty $.
\end{lemma}

Lemma \ref{l1.5} is a consequence of Proposition \ref{p1.2}. To prove Lemma
\ref{l1.6}, which is a generalization of Lemma \ref{l1.4}, it is sufficient
to use Lemma \ref{l1.4}, precompactness of the set $\bigcup\limits_{t\in {\bf
R}}g(t)\subset {\mathcal H}$ for every function $g\in CAP\, ({\bf R},{\mathcal 
H})$, and the Fr\'echet Theorem.

\begin{lemma} \label{l1.6}
Let $f\in B({\bf R},{\mathcal U})$. Then for any $\epsilon ,\delta >0$ there 
are points $x_j\in {\mathcal U}$, $j=1,\dots ,N$ (where $N\in {\bf N}$), such that 
$$
\widetilde {\kappa}\, (\{ t\in {\bf R}:f(t)\in \bigcup\limits_{j=1}^NU_{\delta}
(x_j)\} )<\epsilon \, .
$$
\end{lemma}

\begin{corollary} \label{c1.2}
Let $f\in B({\bf R},{\mathcal U})$. Then there exist points $x_j\in {\mathcal 
U}$, $j\in {\bf N}$, such that 

(1)\ \ ${\rm meas}\, \{ \, t\in {\bf R}:f(t)\notin \overline {\bigcup\limits_{j\in
{\bf N}}x_j}\, \} =0\, $,

(2)\ \ for all $\delta >0$
$$
\widetilde {\kappa}\, (\{ t\in {\bf R}:f(t)\in \bigcup\limits_{j=1}^NU_{\delta}
(x_j)\} )\to 0  \eqno (1.1)
$$
as $N\to +\infty $.
\end{corollary}  

\begin{lemma} \label{l1.7}
Let $f_1\, ,f_2\in B({\bf R},{\mathcal H})$, then $f_1+f_2\in B({\bf R},{\mathcal 
H})$ and ${\rm Mod}\, (f_1+f_2)\subseteq {\rm Mod}\, f_1+{\rm Mod}\, f_2\, $. If 
$f\in B({\bf R},{\mathcal H})$ and $g\in B({\bf R},{\bf C})$, then also $gf\in 
B({\bf R},{\mathcal H})$ and ${\rm Mod}\, gf\subseteq {\rm Mod}\, f+{\rm Mod}\, g$.
\end{lemma}

For $h\in ({\mathcal H},\| .\| )$ we set
$$
{\rm sgn}\, h=\left\{
\begin{array}{ll}
\frac h{\| h\| }\, , &{\text {if}}\ h\neq 0\, ,\\ [0.2cm]
0\, ,   &{\text {if}}\ h=0\, .
\end{array}
\right.
$$

\begin{lemma} \label{l1.8}
Let $f\in B({\bf R},{\mathcal H})$. Suppose that
$$
\widetilde {\kappa}\, (\{ t\in {\bf R}:\| f(t)\| \geq \delta \} )\to 0  \eqno (1.2)
$$
as $\delta \to +0$. Then ${\rm sgn}\, f(.)\in B_1({\bf R},{\mathcal H})$ and
${\rm Mod}\, {\rm sgn}\, f(.)\subseteq {\rm Mod}\, f(.)$ (moreover, for the
set $T=\{ t\in {\bf R}:f(t)=0\} $ we have $\| \chi _T(.)\| ^{(B)}_1=0$).
\end{lemma}

\begin{proof}
For all $j\in {\bf N}$ let us define functions $f_j(t)\doteq j{\mathcal F}
^{\, 1/j}_{\mathcal H}(f(t))$, $t\in {\bf R}$. Lemma \ref{l1.2} implies that
$f_j\in B_1({\bf R},{\mathcal H})$ and ${\rm Mod}\, f_j\subseteq {\rm Mod}\, 
f$. On the other hand, from the condition (1.2) it follows that $\| \chi 
_T(.)\| ^{(B)}_1=0$ and $\| {\rm sgn}\, f(.)-f_j(.)\| ^{(B)}_1\to 0$ as
$j\to +\infty $. Hence ${\rm sgn}\, f(.)\in B_1({\bf R},{\mathcal H})$ and
${\rm Mod}\, {\rm sgn}\, f(.)\subseteq \sum\limits_j{\rm Mod}\, f_j\subseteq
{\rm Mod}\, f(.)$.
\end{proof}

For functions $f,g\in {\mathcal M}_p({\bf R},{\mathcal U})$, $p\geq 1$, we use
the notation
$$
\beta _p(f,g)=\lim\limits_{\delta \to +0}\ \biggl( \ \sup\limits_{T\, \subseteq
\, {\bf R}\, :\, \widetilde {\kappa}\, ({\bf R}\backslash T)\, \leq \, \delta }\
\overline {\lim\limits_{b\to +\infty}}\, \ \frac 1{2b}\ \int\limits_{T\, \cap \,
[-b,b]}\rho ^{\, p}(f(t),g(t))\, dt\, \biggr) ^{1/p}.
$$
The inequality $\beta _p(f,g)\leq D^{(B)}_p(f,g)$ holds for all functions
$f,g\in {\mathcal M}_p({\bf R},{\mathcal U})$. Let
$$
{\mathcal M}^0_p({\bf R},{\mathcal U})\doteq \{ f\in
{\mathcal M}_p({\bf R},{\mathcal U}):\beta _p(f(.),x_0(.))=0\} \, ,
$$
where $x_0(t)\equiv x_0\, $, $t\in {\bf R}$ (the set ${\mathcal M}^0_p({\bf 
R},{\mathcal U})$ does not depend on the choice of the
point $x_0\in {\mathcal U}$); $L^{\infty}({\bf R},{\mathcal U})\subseteq
{\mathcal M}^0_p({\bf R},{\mathcal U})$.

\begin{lemma} \label{l1.9}
For all $p\geq 1$ 
$$
B_p({\bf R},{\mathcal U})=B({\bf R},{\mathcal U})\, \bigcap \, {\mathcal M}
^0_p({\bf R},{\mathcal U})\, .
$$
\end{lemma}

\begin{proof} 
We have $B_p({\bf R},{\mathcal U})\subseteq B({\bf R},{\mathcal U})$. By
Proposition \ref{p1.2}, for any function $f\in B_p({\bf R},{\mathcal U})$
and any $\epsilon >0$ there is a function $f_{\epsilon}\in S_p({\bf R},{\mathcal U})
\cap L^{\infty}({\bf R},{\mathcal U})$ such that $D^{(B)}_p(f,f_{\epsilon})<
\epsilon $. Hence
$$
\beta _p(f(.),x_0(.))\leq \beta _p (f(.),f_{\epsilon}(.))+\beta _p(f_{\epsilon}(.),
x_0(.))\leq
$$ $$
\leq \beta _p (f(.),f_{\epsilon}(.))\leq D^{(B)}_p(f,f_{\epsilon})<\epsilon \, ,
$$
and therefore, $\beta _p(f(.),x_0(.))=0$. The embedding $B_p({\bf R},{\mathcal U})
\subseteq {\mathcal M}^0_p({\bf R},{\mathcal U})$ is proved. Let us now prove
the embedding $B({\bf R},{\mathcal U})\cap {\mathcal M}^0_p({\bf R},{\mathcal U})
\subseteq B_p({\bf R},{\mathcal U})$. By the Fr\'echet Theorem, we can consider
the space ${\mathcal U}=({\mathcal H},\| .\| )$ to be Banach space. Let $f\in
B({\bf R},{\mathcal H})\cap {\mathcal M}^0_p({\bf R},{\mathcal H})\subset
{\mathcal M}_p({\bf R},{\mathcal H})$. The Lemma \ref{l1.4} and the definition
of the set ${\mathcal M}^0_p({\bf R},{\mathcal H})$ imply that for any $\epsilon
>0$ there exists a number $a=a(\epsilon ,f)>0$ such that
$$
\overline {\lim\limits_{b\to +\infty}}\ \ \biggl( \ \frac 1{2b}\ \int\limits
_{ \{ \, t\, \in \, [-b,b]\, :\, \| f(t)\| \, \geq \, a\, \} }\| f(t)\| ^p\, dt
\, \biggr) ^{1/p}<\frac {\epsilon}2\, .
$$
Then ${\mathcal F}_{\mathcal H}^{\, a}(f(.))\in B({\bf R},{\mathcal H})\cap L^{\infty}
({\bf R},{\mathcal H})\subset B_p({\bf R},{\mathcal H})$ and
$$
\| f(.)-{\mathcal F}_{\mathcal H}^{\, a}(f(.))\| ^{(B)}_p\leq
$$ $$
\leq \, \overline {\lim\limits_{b\to +\infty}}\ \ \biggl( \ \frac 1{2b}\ \int\limits
_{ \{ \, t\, \in \, [-b,b]\, :\, \| f(t)\| \, >\, a\, \} }\| f(t)\| ^p\, dt
\, \biggr) ^{1/p}<\frac {\epsilon}2\, .
$$
On the other hand, there is a function $f_{a,\, \epsilon}\in CAP\, ({\bf R},
{\mathcal H})$ such that $\| f_{a,\, \epsilon}\| _{\infty}\leq a$ and $\|
{\mathcal F}_{\mathcal H}^{\, a}(f(.))-f_{a,\, \epsilon}(.)\| ^{(B)}_p< \frac 
{\epsilon}2\, $. Therefore,
$$
\| f-f_{a,\, \epsilon}\| ^{(B)}_p\leq
$$ $$
\leq \| f(.)-{\mathcal F}_{\mathcal H}^{\, a}(f(.))\| ^{(B)}_p+\| {\mathcal F}
_{\mathcal H}^{\, a}(f(.))-f_{a,\, \epsilon}(.)\| ^{(B)}_p< \frac {\epsilon}2 +\frac 
{\epsilon}2=\epsilon \, .
$$
Since the number $\epsilon >0$ can be chosen arbitraryly small, we obtain from
above that $f\in B_p({\bf R},{\mathcal H})$.
\end{proof}

Let $({\rm cl}_b\, {\mathcal U}, {\rm dist})$ be the metric space of
non-empty closed bounded subsets $A\subseteq {\mathcal U}$ with the Hausdorff 
metric
$$
{\rm dist}\, (A,B)={\rm dist}_{\rho}(A,B)=\max \, \{ \sup\limits_{x\in A}
\rho (x,B),\sup\limits_{x\in B}\rho (x,A)\} \, ,\ A,B\in {\rm cl}_b\, {\mathcal 
U}\, ,
$$
where $\rho (x,F)=\inf\limits_{y\in F}\rho (x,y)$ is the distance from a point
$x\in {\mathcal U}$ to a non-empty set $F\subseteq {\mathcal U}$. The metric 
space $({\rm cl}_b\, {\mathcal U},{\rm dist})$ is complete. Let ${\rm cl}\, 
{\mathcal U}$ be the collection of non-empty closed subsets $A\subseteq {\mathcal 
U}$. On the ${\rm cl}\, {\mathcal U}={\rm cl}_b \, ({\mathcal U},\rho ^{\, \prime})$
we define the Hausdorff metric ${\rm dist}_{\rho ^{\, \prime}}$ corresponding to the 
metric $\rho ^{\, \prime}$. The metric space $({\rm cl}\, {\mathcal U},{\rm dist}
_{\rho ^{\, \prime}})$ is also complete. Since  ${\rm dist}^{\, \prime}(A,B)\doteq 
\min \, \{ 1,{\rm dist}\, (A,B)\}={\rm dist}_{\rho ^{\, \prime}}(A,B)$ for all 
$A,B\in {\rm cl}_b \, {\mathcal U}$, it follows that the embedding $({\rm cl}_b\, 
{\mathcal U},{\rm dist}^{\, \prime})\subseteq ({\rm cl}\, {\mathcal U},{\rm dist}
_{\rho ^{\, \prime}})$ is isometric. We define the spaces $B({\bf R},{\rm cl}_b \, 
{\mathcal U})$ and $B_p({\bf R},{\rm cl}_b \, {\mathcal U})$, $p\geq 1$, of {\it 
Besicovitch a.p. multivalued maps} ${\bf R}\ni t\to F(t)\in {\rm cl}_b\, {\mathcal 
U}$ as the spaces of Besicovitch a.p. functions taking values in the metric space 
$({\rm cl}_b \, {\mathcal U},{\rm dist})$. Let $B({\bf R},{\rm cl}\, {\mathcal U})
\doteq B_1({\bf R},({\rm cl}\, {\mathcal U},{\rm dist}_{\rho ^{\, \prime}}))$. The 
following embeddings $B_p({\bf R},{\rm cl}_b \, {\mathcal U})\subseteq B_1({\bf R},
{\rm cl}_b\, {\mathcal U})\subseteq B({\bf R},{\rm cl}_b\, {\mathcal U})\subseteq 
B({\bf R},{\rm cl}\, {\mathcal U})$ hold.

\section{Main results}

Let $B({\bf R})$ be the collection of measurable subsets $T\subseteq {\bf R}$ 
such that $\chi _T\in B_1({\bf R},{\bf R})$. For sets $T\in B({\bf R})$ let 
${\rm Mod}\, T\doteq {\rm Mod}\, \chi _T\, $.

\begin{lemma} \label{l2.1}
Let $T_1\, ,T_2\in B({\bf R})$. Then $T_1\bigcup T_2\in B({\bf R})$, $T_1\bigcap 
T_2\in B({\bf R})$, $T_1\backslash T_2\in B({\bf R})$ and modules ${\rm Mod}\, 
T_1\bigcup T_2\, $, ${\rm Mod}\, T_1\bigcap T_2$ and ${\rm Mod}\, T_1\backslash T_2 
$ are subsets (subgroups) of ${\rm Mod}\, T_1+{\rm Mod}\, T_2\, $.
\end{lemma}

For an arbitrary module $\Lambda \subseteq {\bf R}$ let
${\mathfrak M}^{\, (B)}(\Lambda )$ be the set of sequences $T_j\subseteq {\bf
R}$, $j\in {\bf N}$, of disjoint sets $T_j\in B({\bf R})$ such that ${\rm Mod}\, 
T_j\subseteq \Lambda $, ${\rm meas}\, {\bf R}\backslash \bigcup\limits_jT_j=0$ 
and $\widetilde {\kappa}\, (\bigcup\limits_{j\leq n}T_j)\to 0$ as $n\to 
+\infty $. We shall also assume that the set ${\mathfrak M}^{\, (B)}(\Lambda )$ 
includes the corresponding finite sequences $T_j\, $, $j=1,\dots ,N\, $, which 
can always be supplemented by empty sets to form denumerable ones. The sets 
of sequences $\{ T_j\} \in {\mathfrak M}^{\, (B)}(\Lambda )$ will also be enumerated
by means of several indices.

\begin{lemma} \label{l2.2}
Let  $\Lambda $ be a module in ${\bf R}$ and let $\{ T^{(s)}_j\}\in 
{\mathfrak M}^{\, (B)}(\Lambda )$, $s=1,2$. Then $\{ T^{(1)}_j\bigcap 
T^{(2)}_k\} _{j,\, k}\in {\mathfrak M}^{\, (B)}(\Lambda )$.
\end{lemma}

Lemma \ref{l2.2} is a consequence of Lemma \ref{l2.1}.

Let $\{ T_j\}\in {\mathfrak M}^{\, (B)}(\Lambda )$ and let $J\subseteq {\bf N}$
be an arbitrary non-empty set. Then $\bigcup\limits_{j\in J}T_j\in B({\bf R})$ and
${\rm Mod}\, \bigcup\limits_{j\in J}T_j\subseteq \sum\limits_{j\in J}{\rm Mod}\,
T_j\, $. If $\| \chi _{T_j}\| ^{(B)}_1=0$ for all $j\in J$, then also $\| \chi
_{\bigcup\limits_{j\in J}T_j}\| ^{(B)}_1=0$.

The following Lemma \ref{l2.3} is a consequence of Lemma \ref{l1.7} and the
Fr\'echet Theorem.

\begin{lemma} \label{l2.3}
Suppose that $\{ T_j\} \in {\mathfrak M}^{\, (B)}({\bf R})$ and $f_j\in B({\bf 
R},{\mathcal U})$, $j\in {\bf N}$. Then
$$
\sum\limits_jf_j(.)\chi _{T_j}(.)\in B({\bf R},{\mathcal U})
$$
and
$$
{\rm Mod}\, \sum\limits_jf_j(.)\chi _{T_j}(.)\subseteq \sum\limits_j{\rm Mod}\, f_j
+\sum\limits_j{\rm Mod}\, T_j\, .  \eqno (2.1)
$$
\end{lemma}

\begin{remark} \label{r2.1}
Under the assumptions of Lemma \ref{l2.3}, for indices $j\in {\bf 
N}$ such that $\| \chi _{T_j}\| ^{(B)}_1=0$ (in this case ${\rm Mod}\, T_j=\{ 
0\} $) we can choose arbitrary functions $f_j\in M({\bf R},{\mathcal U})$ and 
delete these indices in the summation on the right-hand side of inclusion (2.1).
\end{remark}

\begin{theorem} \label{th2.1}
Let $f\in B({\bf R},{\mathcal U})$. Then for any $\epsilon >0$ there exist a
sequence $\{ T_j\} _{j\in {\bf N}}\in {\mathfrak M}^{\, (B)}({\rm Mod}\, f)$ and 
points $x_j\in {\mathcal U}$, $j\in {\bf N}$, such that $\rho (f(t),x_j)<\epsilon 
$ for all $t\in T_j\, $, $j\in {\bf N}$.
\end{theorem}

Theorem \ref{th2.1} is proved in Section 3. This Theorem plays a key role in
the paper. Analogous results (on uniform approximation by elementary a.p.
functions) for Stepanov and Weyl a.p. functions were obtained in \cite{2,4}
and \cite{6,7} respectively. For Stepanov a.p. functions stronger assertions
(including a.p. variant of the Lusin Theorem) are contained in \cite{13} and
\cite{14,15} (in last two papers Stepanov a.p. functions are also considered
on relative Bohr compacts).

\begin{corollary} \label{c2.1}
Let $f\in B({\bf R},{\bf R})$. Then for any $a\in {\bf R}$ and $\epsilon >0$ 
there is a set $T\in B({\bf R})$ such that ${\rm Mod}\, T\subseteq {\rm Mod}\, f$, 
$f(t)<a+\epsilon $ for all $t\in T$ and $f(t)>a$ for a.e. $t\in {\bf R}\backslash 
T$. 
\end{corollary}

\begin{theorem} \label{th2.2} 
Let $({\mathcal U},\rho )$ be a complete metric space, let $F\in B({\bf R},{\rm 
cl}\, {\mathcal U})$ and let $g\in B({\bf R},{\mathcal U})$. Then for any 
$\epsilon >0$ there exists a function $f\in B({\bf R},{\mathcal U})$ such that 
${\rm Mod}\, f\subseteq {\rm Mod}\, F+{\rm Mod}\, g\, $, $f(t)\in F(t)$ a.e. and 
$\rho (f(t),g(t))<\rho (g(t),F(t))+\epsilon $ a.e. If, moreover, $F\in B_p({\bf 
R},{\rm cl}_b\, {\mathcal U})$ for some $p\geq 1$, then also $f\in B_p({\bf R},
{\mathcal U})$.
\end{theorem}

\begin{proof} 
Let number $\epsilon \in (0,1]$ be fixed. We choose numbers $\gamma _n>0$, $n\in 
{\bf N}$, such that
$$
\sum\limits_{n=1}^{+\infty}(\gamma _n+\gamma _{n+1})<\frac 16\, .
$$
From Lemmas \ref{l2.1}, \ref{l2.2} and Theorem \ref{th2.1} it follows that for 
each $n\in {\bf N}$ there exist sets 
$F^{(n)}_j\in {\rm cl}\, {\mathcal U}$, points $g^n_j\in {\mathcal U}$ and 
disjoint measurable (in the Lebesgue sense) sets $T^{(n)}_j\subseteq {\bf R}$, 
$j\in {\bf N}$, such that $\{ T^{(n)}_j\} _{j\in {\bf N}}\in {\mathfrak M}^{\, (B)}
({\rm Mod}\, F+{\rm Mod}\, g)$, the functions $F(t)$ and $g(t)$ are defined for
all $t\in \bigcup\limits_jT^{(n)}_j$, and for all $t\in T^{(n)}_j$, $j\in {\bf N}$, 
we have ${\rm dist}_{\rho ^{\, \prime}}(F(t),F^{(n)}_j)<\gamma _n\epsilon <1$ and 
$\rho (g(t),g^n_j)<\gamma _n\epsilon \, $. Let $T=\bigcap\limits_n\bigcup\limits_j
T^{(n)}_j\, $; ${\rm meas}\, {\bf R}\backslash T=0$. By Lemma \ref{l2.2}, for every 
$n\in {\bf N}$
$$
\{ T^{(1)}_{j_1}\bigcap \dots \bigcap T^{(n)}_{j_n}\} _{j_s\, \in \, {\bf N}\, ,\ s=1,
\dots ,n}\in {\mathfrak M}^{\, (B)}({\rm Mod}\, F+{\rm Mod}\, g)\, .
$$
With each number $n\in {\bf N}$ and each collection $\{ j_1\, ,\dots ,\, j_n\} $ 
of indices $j_s\in {\bf N}$, $s=1,\dots ,n$, if $T^{(1)}_{j_1}\bigcap \dots 
\bigcap T^{(n)}_{j_n}\neq \emptyset $, we associate some point $f_{j_1\dots j_n}
\in F^{(n)}_{j_n}\subseteq {\mathcal U}$. These points are determined successively
for $n=1,2,\dots $. For $n=1$ we choose points $f_{j_1}\in F^{(1)}_{j_1}$ such
that the inequalities
$$
\rho (f_{j_1}\, ,g^1_{j_1})<\frac {\epsilon}6+\rho (g^1_{j_1}\, ,F^{(1)}_{j_1})
$$
hold. If points $f_{j_1\dots j_{n-1}}\in F^{(n-1)}_{j_{n-1}}$ have been found for
some $n\geq 2$, we choose points $f_{j_1\dots j_{n-1}j_n}\in F^{(n)}_{j_n}$ 
such that
$$
\rho (f_{j_1\dots j_{n-1}}\, ,f_{j_1\dots j_{n-1}j_n})=\rho ^{\, \prime}(f_{j_1
\dots j_{n-1}}\, ,f_{j_1\dots j_{n-1}j_n})\leq  \eqno (2.2)
$$ $$
\leq 2\, {\rm dist}_{\rho ^{\, \prime}}(F^{(n-1)}_{j_{n-1}},F^{(n)}_{j_n})<2\,
(\gamma _{n-1}+\gamma _n)\epsilon <\frac {\epsilon}3\leq \frac 13\, .
$$
Now let us define functions
$$
f(n;t)=\sum\limits_{j_1\, ,\dots ,\, j_n}f_{j_1\dots j_n}\chi _{T^{(1)}_{j_1}
\bigcap \dots \bigcap T^{(n)}_{j_n}}(t)\, ,\ t\in T\, ,\ n\in {\bf N}\, .
$$
According to Lemmas \ref{l2.2} and \ref{l2.3}, we have $f(n;.)\in B({\bf R},
{\mathcal U})$ and ${\rm Mod}\, f(n;.)\subseteq {\rm Mod}\, F+{\rm Mod}\, g\, $. 
It follows from (2.2) that the inequality
$$
\rho (f(n-1;t),f(n;t))< 2\, (\gamma _{n-1}+\gamma _n)\epsilon \eqno (2.3)
$$
holds for all $t\in T$ and $n\geq 2$. Since the metric space ${\mathcal U}$ is
complete, we obtain from (2.3) that the sequence of functions $f(n;.)$, $n\in {\bf 
N}$, converges as $n\to +\infty $ uniformly on the set $T\subseteq {\bf R}$ 
(therefore, in the metric $D^{\, (B)}$ as well) to a function $f(.)\in B({\bf 
R},{\mathcal U})$ for which ${\rm Mod}\, f\subseteq \sum\limits_n{\rm Mod}\, f(n;.)
\subseteq {\rm Mod}\, F+{\rm Mod}\, g\, $. We have $f(n;t)\in F^{(n)}_{j_n}$ and ${\rm 
dist}_{\rho ^{\, \prime}}(F(t),F^{(n)}_{j_n})<\gamma _n\epsilon <\frac 16$ for 
all $t\in T^{(n)}_{j_n}\bigcap T$. Since $\gamma _n\to 0$ as $n\to +\infty $, it
follows from this that $f(t)\in F(t)$ for all $t\in T$ (for a.e. $t\in {\bf R}$). 
With each number $t\in T$ we associate an infinite collection of indices $\{ j_1
\, ,\dots ,\, j_n\, ,\dots \} $ in such a way that $t\in T^{(n)}_{j_n}$, $n\in 
{\bf N}$. Then (for all $t\in T$)
$$
\rho (f(t),g(t))\leq \sum\limits_{n=1}^{+\infty}\rho (f_{j_1\dots j_n}\, ,f_{j_1
\dots j_nj_{n+1}})+\rho (f_{j_1}\, ,g^1_{j_1})+\rho (g^1_{j_1}\, ,g(t))<
$$ $$
<2\, \sum\limits_{n=1}^{+\infty}(\gamma _n+\gamma _{n+1})\epsilon +\frac 
{\epsilon}3+\rho (g^1_{j_1}\, ,F^{(1)}_{j_1})<
$$ $$
<\frac {2\epsilon }3+|\rho (g^1_{j_1}\, ,F^{(1)}_{j_1})-\rho (g^1_{j_1}\, 
,F(t))|+|\rho (g^1_{j_1}\, ,F(t))-\rho (g(t),F(t))|+\rho (g(t),F(t))<
$$ $$
<\frac {2\epsilon }3+\gamma _1\epsilon +\gamma _1\epsilon +\rho (g(t),F(t))<
\epsilon +\rho (g(t),F(t))\, .
$$
If $F\in B_p({\bf R},{\rm cl}_b\, {\mathcal U})\subseteq B({\bf R},{\rm cl}\, 
{\mathcal U})$, $p\geq 1$, then $f\in B_p({\bf R},{\mathcal U})$. Indeed, for 
a.e. $t\in {\bf R}$ we have
$$
\rho (f(t),x_0)\leq \sup\limits_{x\, \in \, F(t)}\rho (x,x_0)={\rm dist}\, (F(t),
\{ x_0\} )\, ,
$$
furthermore, ${\rm dist}\, (F(.),\{ x_0\} )\in {\mathcal M}^0_p({\bf R},{\bf R})$. 
Hence (see Lemma \ref{l1.9}) $f(.)\in B({\bf R},{\mathcal U})\bigcap {\mathcal 
M}^0_p({\bf R},{\mathcal U})=B_p({\bf R},{\mathcal U})$.
\end{proof}

\begin{corollary} \label{c2.2}
Let $({\mathcal U},\rho )$ be a complete separable metric space and let $F\in 
B({\bf R},{\rm cl}\, {\mathcal U})$. Then there exist functions $f_j\in B({\bf R},
{\mathcal U})$, $j\in {\bf N}$, such that ${\rm Mod}\, f_j\subseteq {\rm Mod}\, F$ 
and $F(t)=\overline {\bigcup\limits_jf_j(t)}$ for a.e. $t\in {\bf R}$ (if $F\in 
B_p({\bf R},{\rm cl}_b\, {\mathcal U})\subseteq B({\bf R},{\rm cl}\, {\mathcal U})$, 
$p\geq 1$, then all functions $f_j$ belong to the space $B_p({\bf R},{\mathcal U})$).
\end{corollary}

\begin{proof} 
Let us choose points $x_k\in {\mathcal U}$, $k\in {\bf N}$, which form a 
countable dense set of the metric space ${\mathcal U}$. By Theorem \ref{th2.2},
for all $k,n\in {\bf N}$ there are functions $f_{k,\, n}\in B({\bf R},{\mathcal 
U})$ such that ${\rm Mod}\, f_{k,\, n}\subseteq {\rm Mod}\, F$, $f_{k,\, n}(t)\in F(t)$ 
a.e. and $\rho (f_{k,\, n}(t),x_k)<2^{-n}+\rho (x_k\, ,F(t))$ a.e. Furthermore, in
the case $F\in B_p({\bf R},{\rm cl}_b\, {\mathcal U})$, $p\geq 1$, we also have
$f_{k,\, n}\in B_p({\bf R},{\mathcal U})$, $k,n\in {\bf N}$. It remains to
renumber the functions $f_{k,\, n}$ by a single index $j\in {\bf N}$.
\end{proof}

The proof of following Theorem \ref{th2.3} is analogous to the proof of Theorem
1.3 in \cite{6} (in which Weyl a.p. functions and multivalued maps were
considered). To prove Theorem \ref{th2.3} it is necessary to use Theorem
\ref{th2.2}, Corollary \ref{c2.1} and Lemmas \ref{l2.1} and \ref{l2.3}. Analogous
(to Theorem \ref{th2.3}) result for Stepanov a.p. functions and multivalued maps
can be found in \cite{14}.

\begin{theorem} \label{th2.3}
Let $({\mathcal U},\rho )$ be a complete metric space, let $F\in B({\bf R},{\rm 
cl}\, {\mathcal U})$ and let $g\in B({\bf R},{\mathcal U})$. Then for any
non-decreasing function $[0,+\infty )\ni t\to \eta (t)\in {\bf R}$, for which $\eta 
(0)=0$ and $\eta (t)>0$ for all $t>0$, there exists a function $f\in B({\bf R},
{\mathcal U})$ such that ${\rm Mod}\, f \subseteq {\rm Mod}\, F+{\rm Mod}\, g$, 
$f(t)\in F(t)$ a.e. and $\rho (f(t),g(t))\leq \rho (g(t),F(t))+\eta (\rho (g(t),
F(t)))$ a.e. Moreover, if $F\in B_p({\bf R},{\rm cl}_b\, {\mathcal U})\subseteq 
B({\bf R},{\rm cl}\, {\mathcal U})$, $p\geq 1$, then $f\in B_p({\bf R},{\mathcal 
U})$.
\end{theorem}

The following Theorems can be also proved (using Theorems \ref{th2.1}, \ref{th2.2}
and Lemmas \ref{l2.1}, \ref{l2.2} and \ref{l2.3}) by analogy with appropriate
assertions on Stepanov \cite{5,14} and Weyl \cite{7} a.p. functions and
multivalued maps.

The points $x_j\in {\mathcal U}$, $j=1,\dots ,n$, are said to form {\it
$\epsilon$-net} for (non-empty) set $F\subseteq {\mathcal U}$, $\epsilon >0$, if
$F\subseteq \bigcup\limits_jU_{\epsilon}(x_j)$.

\begin{theorem} \label{th2.4}
Let $({\mathcal U},\rho )$ be a complete metric space, let $F\in B({\bf R},{\rm 
cl}_b\, {\mathcal U})$ and let $\epsilon >0$, $n\in {\bf N}$. Suppose that for
a.e. $t\in {\bf R}$ there are points $x_j(t)\in F(t)$, $j=1,\dots ,n$, which
form $\epsilon$-net for the set $F(t)$. Then for any $\epsilon ^{\, \prime}>
\epsilon $ there exist functions $f_j\in B({\bf R},{\mathcal U})$, $j=1,\dots 
,n$, such that ${\rm Mod}\, f_j\subseteq {\rm Mod}\, F$, $f_j(t)\in F(t)$ a.e.
and for a.e. $t\in {\bf R}$ the points $f_j(t)$, $j=1,\dots ,n$, form $\epsilon
^{\, \prime}$-net for the set $F(t)$.
\end{theorem}

\begin{corollary} \label{c2.3}
Let $({\mathcal U},\rho )$ be a compact metric space. Then a multivalued map
${\bf R}\ni t\to F(t)\in {\rm cl}\, {\mathcal U}={\rm cl}_b\, {\mathcal U}$
belongs to the space $B({\bf R},{\rm cl}\, {\mathcal U})=B_1({\bf R},{\rm 
cl}_b\, {\mathcal U})$ if and only if for each $\epsilon >0$ there exist a
number $n\in {\bf N}$ and functions $f_j\in B({\bf R},{\mathcal U})=B_1({\bf R},
{\mathcal U})$, $j=1,\dots ,n$, such that $f_j(t)\in F(t)$ a.e. and points 
$f_j(t)$, $j=1,\dots ,n$, for a.e. $t\in {\bf R}$ form $\epsilon$-net for the set 
$F(t)$ (furthermore, the functions $f_j$ for the multivalued map $F\in B({\bf R},
{\rm cl}\, {\mathcal U})$ can be chosen in such a way that ${\rm Mod}\, f_j
\subseteq {\rm Mod}\, F$).
\end{corollary}

\begin{theorem} \label{th2.5}
Let $({\mathcal U},\rho )$ be a compact metric space. Then a multivalued map
${\bf R}\ni t\to F(t)\in {\rm cl}\, {\mathcal U}$ belongs to the space $B({\bf R},
{\rm cl}\, {\mathcal U})$ if and only if there exist functions $f_j\in B({\bf R},
{\mathcal U})$, $j\in {\bf N}$, such that $F(t)=\overline {\bigcup\limits_jf(t)}$ 
a.e. and the set $\{ f_j(.): j\in {\bf N}\} $ is precompact in the metric space
$L^{\infty}({\bf R},{\mathcal U})$ (furthermore, the functions $f_j$ for the 
multivalued map $F\in B({\bf R},{\rm cl}\, {\mathcal U})$ can be chosen in such 
a way that ${\rm Mod}\, f_j\subseteq {\rm Mod}\, F$).
\end{theorem}

For non-empty set $F\subseteq {\mathcal U}$ we shall use the notation
$F^{\, \delta}=\{ x\in {\mathcal U}:\rho (x,F)<\delta \} $, $\delta >0$.

\begin{theorem} \label{th2.6}
Let $({\mathcal U},\rho )$ be a complete metric space, let $F\in B({\bf R},{\rm 
cl}_b\, {\mathcal U})$, $\epsilon >0$, $\delta >0$, $n\in {\bf N}$, and let 
$g_j\in B({\bf R},{\mathcal U})$, $j=1,\dots ,n$. Suppose that for
a.e. $t\in {\bf R}$ the set of points $x_j(t)=g_j(t)$, for which $g_j(t)\in
(F(t))^{\delta}$, can be supplemented (if it consists of less than $n$ points)
to $n$ points $x_j(t)\in (F(t))^{\delta}$, $j=1,\dots ,n$, which
form $\epsilon$-net for the set $F(t)$ (coincident points with different
indices are considered here as different points). Then for any $\epsilon ^{\, \prime}>
\epsilon +\delta $ there exist functions $f_j\in B({\bf R},{\mathcal U})$, $j=1,
\dots ,n$, such that ${\rm Mod}\, f_j\subseteq {\rm Mod}\, F+\sum\limits_{k=1}^n
{\rm Mod}\, g_k\, $, $f_j(t)\in F(t)$ a.e., $f_j(t)=g_j(t)$ for a.e. $t\in
\{ \tau \in {\bf R}:g_j(\tau )\in F(\tau )\} $ and the points $f_j(t)$, $j=1,
\dots ,n$, for a.e. $t\in {\bf R}$ form $\epsilon^{\, \prime}$-net for the set $F(t)$.
\end{theorem}

Let $({\mathcal U},\rho )$ and $({\mathcal V},\rho _{\mathcal V})$ be complete
metric spaces and let $C({\mathcal U},{\mathcal V})$ be the space of continuous
functions ${\mathcal F}:{\mathcal U}\to {\mathcal V}$ endowed with metric
$$
d_{C({\mathcal U},{\mathcal V})}({\mathcal F}_1\, ,{\mathcal F}_2)=\sup\limits
_{x\, \in \, {\mathcal U}}\, \min \, \{ 1,\rho _{\mathcal V}({\mathcal F}_1(x),
{\mathcal F}_2(x))\} \, ,\ {\mathcal F}_1\, ,\, {\mathcal F}_2\in C({\mathcal U},
{\mathcal V})\, .
$$
We denote by ${\mathcal F}(.|_Y)$ the restriction of a function ${\mathcal F}:
{\mathcal U}\to {\mathcal V}$ to a non-empty set $Y\subseteq {\mathcal U}$. In 
following Lemmas we consider the superposition of Besicovitch a.p. functions.

\begin{lemma} \label{l2.4}
Let $({\mathcal U},\rho )$ and $({\mathcal V},\rho _{\mathcal V})$ be complete
metric spaces, let ${\mathcal F}\in C({\mathcal U},{\mathcal V})$ and let $f\in 
B({\bf R},{\mathcal U})$. Then ${\mathcal F}(f(.))\in B({\bf R},{\mathcal V})$ 
and ${\rm Mod}\, {\mathcal F}(f(.))\subseteq {\rm Mod}\, f(.)$.
\end{lemma}

\begin{proof}
We have ${\mathcal F}(f(.))\in M({\bf R},{\mathcal V})={\mathcal M}_1({\bf R},
({\mathcal V},\rho ^{\, \prime}_{\mathcal V}))$. Let $\epsilon \in (0,1]$, $\delta 
>0$. By Theorem \ref{th2.1}, for every $k\in {\bf N}$ there are sequences $\{ 
T^{(k)}_j\} _{j\in {\bf N}}\in {\mathfrak M}^{\, (B)}({\rm Mod}\, f)$ and points 
$x^{(k)}_j\in {\mathcal U}$, $j\in {\bf N}$, such that $\rho (f(t),x^{(k)}_j)<
\frac 1k$ for all $t\in T^{(k)}_j$, $j\in {\bf N}$. Let us choose numbers $j(k)
\in {\bf N}$, $k\in {\bf N}$, for which
$$
\widetilde {\kappa}\, \biggl( \ \bigcup\limits_{j=1}^{j(k)}T^{(k)}_j\, \biggr) <
2^{-k}\epsilon \, .
$$
Let $X_1=\bigcup\limits_{j\, \leq \, j(1)}x^{(1)}_j$. For every $k\in {\bf N}
\backslash \{ 1\} $ we denote by $X_k$ the set of points $x^{(k)}_j$, $j=1,\dots 
,j(k)$, for which for any $k^{\, \prime}=1,\dots ,\, k-1$ there exists a point 
$x^{(k^{\, \prime})}_{j^{\, \prime}}$, $j^{\, \prime}=1,\dots ,j(k^{\, \prime})$, 
such that $\rho (x^{(k)}_j,x^{(k^{\, \prime})}_{j^{\, \prime}})<\frac 1k+\frac 
1{k^{\, \prime}}<\frac 2{k^{\, \prime}}\, $. If $T^{(1)}_{j_1}\bigcap \dots 
\bigcap T^{(k)}_{j_k}\neq \emptyset $, $k\in {\bf N}$, where $j_s\in \{ 1,\dots 
,j(s)\} $, $s=1,\dots ,k$, then $x^{(k)}_{j_k}\in X_k\, $, hence from the
precompactness of the set $\bigcup\limits_{k\, \in \, {\bf N}}X_k\subseteq 
{\mathcal U}$ and from the continuity of the function ${\mathcal F}$ it follows 
that there is a number $k_0\in {\bf N}$ such that for all $j_k=1,\dots ,j(k)$, 
where $k=1,\dots ,k_0\, $, and for all $t,t^{\, \prime}\in T^{(1)}_{j_1}\bigcap 
\dots \bigcap T^{(k_0)}_{j_{k_0}}$ the inequality
$$
\rho _{\mathcal V}({\mathcal F}(f(t)),{\mathcal F}(f(t^{\, \prime})))<\delta 
$$
holds. If $T^{(1)}_{j_1}\bigcap \dots \bigcap T^{(k_0)}_{j_{k_0}}\neq \emptyset $, 
where $j_k\in \{ 1,\dots ,j(k)\} $, $k=1,\dots ,k_0\, $, we choose some 
numbers $t_{j_1\dots j_{k_0}}\in T^{(1)}_{j_1}\bigcap \dots \bigcap T^{(k_0)}_{j
_{k_0}}$. Let
$$
T(k_0)=\bigcap\limits_{k\, =\, 1,\dots ,k_0}\, \bigcup\limits_{j_k=1}^{j(k)}T^{(k)}
_{j_k}\, .
$$
By Lemmas \ref{l2.1} and \ref{l2.3},
$$
{\mathcal G}_{k_0}(.)\doteq \sum\limits_{j_k\, =\, 1,\dots ,j(k);\ k\, =\, 1,
\dots ,k_0}{\mathcal F}(f(t_{j_1\dots j_{k_0}}))\, \chi _{T^{(1)}_{j_1}\bigcap \dots 
\bigcap T^{(k_0)}_{j_{k_0}}}(.)\, +
$$ $$
+\, y_0\, \chi _{{\bf R}\, \backslash T(k_0)}(.)\in B
({\bf R},{\mathcal V})\, ,
$$ 
where $y_0\in {\mathcal V}$, and
$$
{\rm Mod}\, {\mathcal G}_{k_0}(.)\subseteq \sum\limits_{j_k\, =\, 1,\dots ,j(k);\ 
k\, =\, 1,\dots ,k_0}{\rm Mod}\, T^{(k)}_{j_k}\subseteq {\rm Mod}\, f(.)\, .
$$
Furthermore, $\rho _{\mathcal V}({\mathcal F}(f(t)),{\mathcal G}_{k_0}(t))<
\delta $ for all $t\in T(k_0)$ and
$$
\widetilde {\kappa}\, (T(k_0))\leq \sum\limits_{k=1,\dots ,k_0}\, \widetilde 
{\kappa}\, \biggl( \, \bigcup\limits_{j_k=1}^{j(k)}T^{(k)}_{j_k}\, \biggr)
<\sum\limits_{k=1,\dots ,k_0}2^{-k}\epsilon <\epsilon \, .
$$
Hence $D^{\, (B)}({\mathcal F}(f(.)),{\mathcal G}_{k_0}(.))<\epsilon +\delta $. 
Since the numbers $\epsilon >0$ and $\delta >0$ can be chosen arbitraryly small, 
it follows that ${\mathcal F}(f(.))\in B({\bf R},{\mathcal V})$ and ${\rm Mod}\, 
{\mathcal F}(f(.))\subseteq {\rm Mod}\, f(.)$. 
\end{proof}

\begin{lemma} \label{l2.5} 
Let $({\mathcal U},\rho )$ and $({\mathcal V},\rho _{\mathcal V})$ be complete
metric spaces. Suppose that a function ${\bf R}\ni t\to {\mathcal F}(.;t)\in 
C({\mathcal U},{\mathcal V})$ belongs to the space $B_1({\bf R},(C({\mathcal U},
{\mathcal V}),d_{C({\mathcal U},{\mathcal V})}))$ and $f\in B({\bf R},{\mathcal 
U})$. Then ${\mathcal F}(f(.);.)\in B({\bf R},{\mathcal V})$ and ${\rm Mod}\, 
{\mathcal F}(f(.);.)\subseteq {\rm Mod}\, {\mathcal F}(.;.)+{\rm Mod}\, f(.)$.
\end{lemma}

\begin{proof}
Theorem \ref{th2.1} implies that for any $\epsilon >0$ there are a sequence $\{ 
T_j\} _{j\in {\bf N}}\in {\mathfrak M}^{\, (B)}({\rm Mod}\, {\mathcal F}(.;.))$ 
and functions ${\mathcal F}_j\in C({\mathcal U},{\mathcal V})$, $j\in {\bf N}$, 
such that $d_{C({\mathcal U},{\mathcal V})}({\mathcal F}(.;t),{\mathcal F}_j(.))
<\epsilon $ for all $t\in T_j\, $, $j\in {\bf N}$. By Lemmas \ref{l2.3} and 
\ref{l2.4},
$$
\sum\limits_{j\in {\bf N}}\, {\mathcal F}_j(f(.))\, \chi _{T_j}(.)\in B({\bf R},
{\mathcal V})\, ,
$$ $$
{\rm Mod}\, \sum\limits_{j\in {\bf N}}\, {\mathcal F}_j(f(.))\, \chi _{T_j}(.)\subseteq
{\rm Mod}\, {\mathcal F}(.;.)+{\rm Mod}\, f(.)\, .
$$
On the other hand, ${\mathcal F}(f(.);.)\in M({\bf R},{\mathcal V})={\mathcal M}_1
({\bf R},({\mathcal V},\rho ^{\, \prime}_{\mathcal V}))$ and
$$
D^{\, (B)}\bigl( \, {\mathcal F}(f(.);.),\, \sum\limits_{j\in 
{\bf N}}\, {\mathcal F}_j(f(.))\, \chi _{T_j}(.)\, \bigr) <\epsilon \, .
$$
Hence (since the number $\epsilon >0$ can be chosen arbitraryly small), 
${\mathcal F}(f(.);.)\in B({\bf R},{\mathcal V})$ and ${\rm Mod}\, {\mathcal F}
(f(.);.)\subseteq {\rm Mod}\, {\mathcal F}(.;.)+{\rm Mod}\, f(.)$. 
\end{proof}

\begin{remark} \label {r2.2}
From Lemmas \ref{l1.9}, \ref{l2.3}, \ref{l2.4} and Theorem \ref{th2.1}
we obtain also the following assertion. Let $({\mathcal U},\rho )$ and 
$({\mathcal V},\rho _{\mathcal V})$ be complete metric spaces, let $r>0$ and let 
$p\geq 1$. Suppose that a function ${\bf R}\ni t\to {\mathcal F}(.;t)\in 
C({\mathcal U},{\mathcal V})$ satisfies the following two conditions:

(1) for every $x\in {\mathcal U}$ the function ${\bf R}\ni t\to {\mathcal F}
(.|_{U_r(x)};t)\in C(U_r(x),{\mathcal V})$ belongs to the space 
$$
B_1({\bf R},(C(U_r(x),{\mathcal V}),d_{C(U_r(x),{\mathcal V})}))\, ;
$$

(2) there exist a number $C\geq 0$ and a function $C(.)\in {\mathcal M}^0_p
({\bf R},{\bf R})$ such that for a.e.  $t\in {\bf R}$ the inequality
$$
\rho _{\mathcal V}({\mathcal F}(x;t),y_0)\leq A\rho (x,x_0)+B(t)
$$
holds for all $x\in {\mathcal U}$, where $x_0\in {\mathcal U}$ and $y_0\in 
{\mathcal V}$ are some fixed points.

Then for any function $f\in B_p({\bf R},{\mathcal U})$ we have ${\mathcal F}
(f(.);.)\in B_p({\bf R},{\mathcal V})$ and 
$$
{\rm Mod}\, {\mathcal F}(f(.);.)\subseteq {\rm Mod}\, f(.)+\sum\limits_{x\, \in 
\, {\mathcal U}}{\rm Mod}\, {\mathcal F}(.|_{U_r(x)};.)\, .
$$
\end{remark}

\section{Proof of Theorem \ref{th2.1}}

Let ${\mathcal A}^{(B)}$ be the collection of sets $\mathbb F\subset B({\bf R},
{\bf R})$ such that 
$$
\lim\limits_{\tau _0\to +0}\ \, \sup\limits_{f\in {\mathbb F}}\ 
\sup\limits_{\tau \in [0,\tau _0]}\, D^{\, (B)}(f(.),f(.+\tau ))=0\, .
$$
If $f\in B_p({\bf R},{\mathcal U})$, $p\geq 1$, then $D^{\, (B)}_p(f(.),f(.+\tau ))
\to 0$ as $\tau \to 0$, therefore (in particular) for all functions $f\in B({\bf R},
{\bf R})$ we have $\{ f(.)+a:a\in {\bf R}\} \in {\mathcal A}^{(B)}$.

For a measurable set $T\subseteq {\bf R}$ let denote
$$
\kappa \, (T)\doteq \widetilde {\kappa}\, ({\bf R}\backslash T)=\overline 
{\lim\limits_{b\to +\infty}}\ \ \frac 1{2b}\ \, {\rm meas}\, [-b,b]\bigcap T\, .
$$
If $T_1\, ,\, T_2\subseteq {\bf R}$ are measurable sets, then $\kappa \, (T_1\cup
T_2)\leq \kappa \, (T_1)+\kappa \, (T_2)$.

The following Theorem \ref{th3.1} is proved in Section 4 and its special case 
for the set ${\mathbb F}=\{ f\} $, $f\in B({\bf R},{\bf R})$, is essentially 
used in the proof of Theorem \ref{th2.1}.

\begin{theorem} \label{th3.1}
Let ${\mathbb F}\in {\mathcal A}^{(B)}$, $\Delta >0$, $b>0$. Then there exists 
$b$-periodic function $g\in C({\bf R},{\bf R})$, dependent on ${\mathbb F}$, 
$\Delta $ and $b$, for which $\| g\| _{\infty}<\Delta $, such that for every
$\epsilon \in (0,1]$ there is a number $\delta = \delta (\epsilon ,\Delta )>0$
such that for all functions $f\in {\mathbb F}$ the inequality
$$
\kappa \, (\{ t\in {\bf R}:|f(t)+g(t)|<\delta \} )<\epsilon
$$
holds.
\end{theorem}

{\it Proof of Theorem \ref{th2.1}}. If ${\rm Mod}\, f=\{ 0\} $, then there is a 
constant function $f_0(t)\equiv f_0\in {\mathcal U}$, $t\in {\bf R}$, for which 
$D^{\, (B)}(f(.),f_0(.))=0$, therefore, there is a set $T\in B({\bf R})$ 
such that $\widetilde {\kappa}\, (T)=0$ and $\rho (f(t),f_0)<\epsilon $ for all 
$t\in T$. From this (using the measurability of function $f(.)$) we obtain the
assertion to be proved (furthermore, $T_1=T$). Next, suppose that ${\rm Mod}\, f
\neq \{ 0\} $. Let $x_j\in {\mathcal U}$, $j\in {\bf N}$, be the points 
determined in Corollary \ref{c1.2} for the function $f\in B({\bf R},{\mathcal U})$.
By Corollary \ref{c1.1}, for all $j\in {\bf N}$ we have $\rho (f(.),x_j)\in B({\bf 
R},{\bf R})$ and ${\rm Mod}\, \rho (f(.),x_j)\subseteq {\rm Mod}\, f(.)$. We choose
a number $b>0$ such that $\frac {2\pi }b\in {\rm Mod}\, f$. Theorem \ref{th3.1}
implies the existence of $b$-periodic function $g_j(.)\in C({\bf R},{\bf R})$, 
$j\in {\bf N}$, such that $\| g_j\| _{\infty}<\frac {\epsilon }3$ and 
$$
\kappa \, (\{ t\in {\bf R}:|\, \rho (f(t),x_j)-\frac {2\epsilon }3+g_j(t)\, |<\delta 
\} )\to 0
$$
as $\delta \to +0$ (instead of functions $g_j$ we could choose one function $g_0=
g_j\, $, $j\in {\bf N}$, but it doesn't change the proof). Let $T^{\, \prime}_j=\{ 
\, t\in {\bf R}:\rho (f(t),x_j)+g_j(t)\leq \frac {2\epsilon }3\, \} $, $j\in 
{\bf N}$. According to Lemma \ref{l1.8}, we get $T^{\, \prime}_j\in B({\bf R})$ 
and ${\rm Mod}\, T^{\, \prime}_j\subseteq {\rm Mod}\, \rho (f(.),x_j)+\frac {2\pi }b
\, {\bf Z}\subseteq {\rm Mod}\, f(.)$. If $t\in T^{\, \prime}_j\, $, then $\rho 
(f(t),x_j)<\epsilon $. We denote $T_1=T^{\, \prime}_1$ and $T_j=T^{\, \prime}_j
\backslash \bigcup\limits_{k\, <\, j}T^{\, \prime}_k$ for $j\geq 2$. The sets 
$T_j$, $j\in {\bf N}$, are disjoint and $\bigcup\limits_{j\leq N}T_j=\bigcup\limits
_{j\leq N}T^{\, \prime}_j$ for all $N\in {\bf N}$. It follows from Lemma \ref{l2.1} 
that $T_j\in B({\bf R})$, ${\rm Mod}\, T_j\subseteq {\rm Mod}\, f$. Furthermore, 
$\rho (f(t),x_j)<\epsilon $ for all $t\in T_j\, $, $j\in {\bf N}$, and for every 
$N\in {\bf N}$ and a.e. $t\in {\bf R}\backslash \bigcup\limits_{j\leq N}T_j$ we 
have $\rho (f(t),x_j)> \frac {\epsilon}3$ for all $j=1,\dots ,N$. Hence (see 
Corollary \ref{c1.2}) ${\rm meas}\, {\bf R}\, \backslash \bigcup\limits_{j}T_j=0$ 
and (see (1.1) for $\delta =\frac {\epsilon}3\, $) $\widetilde {\kappa}\, \bigl( \,
\bigcup\limits_{j\leq N}T_j\, \bigr) \to 0$ as $N\to +\infty $, that is, 
$\{ T_j\} \in {\mathfrak M}^{\, (B)}({\rm Mod}\, f)$. \hfill $\square$

\section{Proof of Theorem \ref{th3.1}}

\begin{lemma} \label{l4.1}
Let ${\mathbb F}\in {\mathcal A}^{(B)}$, $\Delta >0$. Then for any $\epsilon \in 
(0,1]$ there exist numbers $\delta =\delta (\epsilon ,\Delta )>0$ and $\widetilde 
{\alpha}=\widetilde {\alpha}(\epsilon ,\Delta ,{\mathbb F})>0$ such that for all 
$\alpha \geq \widetilde {\alpha}$ and all functions $f\in {\mathbb F}$
$$
\kappa \, (\{ t\in {\bf R}:|\, f(t)+\Delta \sin \alpha t\, |<\delta \} )< \epsilon \, .
$$
\end{lemma}

\begin{proof}
Let us choose a number $N=N(\epsilon )\in {\bf N}$ for which $(N+1)^{-1}<\frac 
{\epsilon}2$ (then $N\geq 2$). Let 
$$
\epsilon ^{\, \prime}\doteq \frac 12\, \epsilon N^{-1}(N+1)^{-1}\leq \frac 
{\epsilon}{12}<1\, ,\ \delta ^{\, \prime}\doteq 2\sin \frac {\pi}{2N}\, \sin \frac 
{\pi \epsilon ^{\, \prime}}2\, , 
$$ $$
\delta =\delta (\epsilon ,\Delta )\doteq \min \, \{ 1,\frac 13\, \delta ^{\, \prime}
\Delta \} \, . 
$$
There is a number $\tau _0=\tau _0(\epsilon ,\Delta ,{\mathbb F})> 0$ such that the
inequality 
$$
D^{\, (B)}(f(.),f(.+\tau ))<\epsilon ^{\, \prime}\delta 
$$
holds for all $f\in {\mathbb F}$ and $\tau \in [0,\tau _0]$. We define the number
$\widetilde {\alpha}=\pi \tau _0^{-1}$. Let $0<\tau \leq \tau_0\, $, $\alpha 
\doteq \pi \tau ^{-1}\geq \widetilde {\alpha}$. For all $j=1,\dots ,N$ (and all
functions $f\in {\mathbb F}$) let us define the sets
$$
{\mathcal L}_j(f,\tau )=\{ t\in {\bf R}:|\, f(t+\frac jN\, \tau )-f(t)\, |\geq 
\delta \} \, .
$$
We have
$$
\kappa \,({\mathcal L}_j(f,\tau ))\leq \, \frac 1{\delta}\ \, \overline {\lim\limits
_{b\to +\infty}}\ \ \frac 1{2b}\ \int\limits_{-b}^b\min \, \{ 1,|\, f(t+\frac jN\, 
\tau )-f(t)\, |\} \, dt\, =
$$ $$
=\, \frac 1{\delta}\ D^{\, (B)}(f(.),f(.+\frac jN\, \tau ))<\epsilon ^{\, \prime}\, .
$$
For $j=1,\dots ,N$ we also consider the sets  
$$
{\mathcal N}_j(\tau )=\{ t\in {\bf R}:|\, \cos \alpha (t+\frac j{2N}\, \tau )\, 
\sin \frac {\alpha j}{2N}\, \tau \, |\leq \frac {\delta ^{\, \prime}}2\, \} \, .
$$ 
If $t\in {\mathcal N}_j(\tau )$, then
$$
|\, \cos (\alpha t+\frac {j\pi }{2N}\, )\, |\leq \frac 12\ \delta ^{\, \prime} 
\sin ^{-1}\frac {\pi}{2N}=\sin \frac {\pi \epsilon ^{\, \prime}}2\, ,
$$
therefore, the number $t$ belongs to one of closed intervals $[\beta ^-_s\, ,
\beta ^+_s]$, $s\in {\bf Z}$, where
$$
\beta ^{\pm}_s=\bigl( s+\frac 12\, \bigr) \, \frac {\pi}{\alpha}-\frac {j\pi }{2N
\alpha }\pm \frac {\pi \epsilon ^{\, \prime}}{2\alpha }\, ,
$$
and hence $\kappa \, ({\mathcal N}_j(\tau ))\leq \epsilon ^{\, \prime}$. In what 
follows, we shall suppose that the sets ${\mathcal L}_j(f,\tau )$ contain (in
addition) the numbers $t\in {\bf R}$ for which at the least one of the
functions $f(t)$, $f(t+\frac jN\, \tau )$ is not defined (these numbers form the
set of measure zero). Let
$$
{\mathcal L}(f,\tau )=\bigcup\limits_{j=1}^N\, \biggl( \, \bigcup\limits_{s=0}^{N-j}
({\mathcal L}_j(f,\tau )-\frac sN\, \tau )\biggr)
$$
(here ${\mathcal L}_j(f,\tau )-\frac sN\, \tau =\{ t=\eta -\frac sN\, \tau :\eta \in
{\mathcal L}_j(f,\tau )\} $). Since $\kappa \, ({\mathcal L}_j(f,\tau ))<\epsilon ^{\, 
\prime}$, $j=1,\dots ,N$, we get 
$$
\kappa \, ({\mathcal L}(f,\tau ))<\frac 12 \, N(N+1)\epsilon ^{\, \prime}=
\frac {\epsilon}4\, .
$$
If $t\in {\bf R}\backslash {\mathcal L}(f,\tau )$, then 
$$
|\, f(t+\frac {j_1}N\, \tau )-f(t+\frac {j_2}N\, \tau )\, |<\delta 
$$
for all $j_1\, , j_2\in \{ 0,1,\dots ,N\} $. Let
$$
{\mathcal N}(\tau )=\bigcup\limits_{j=1}^N\, \biggl( \, \bigcup\limits_{s=0}^{N-j}
\, (\, {\mathcal N}_j(\tau )-\frac sN\, \tau )\biggr) \, .
$$
Since $\kappa \, ({\mathcal N}_j(\tau ))<\epsilon ^{\, \prime}$, $j=1,\dots ,N$, we 
also get 
$$
\kappa \, ({\mathcal N}(\tau ))<\frac 12 \, N(N+1)\epsilon ^{\, \prime}=
\frac {\epsilon}4\, .
$$
If $t\in {\bf R}\backslash {\mathcal N}(\tau )$, then for all $j_1\, , j_2\in 
\{ 0,1,\dots ,N\} $, $j_1<j_2\, $, we have
$$
|\, \Delta \sin \alpha (t+\frac {j_1}N\, \tau )-\Delta \sin \alpha (t+\frac {j_2}N\, 
\tau )\, |=
$$ $$
=2\Delta \, |\, \cos \alpha (t+\frac {j_1}N\, \tau +\frac {j_2-j_1}{2N}\, \tau )\, \sin
\alpha \, \frac {j_2-j_1}{2N}\, \tau \, |>\Delta \delta ^{\, \prime}\geq 3\delta \, .
$$
Let $G(t)=f(t)+\Delta \sin \alpha t$, $t\in {\bf R}$. We define the set
$$
{\mathcal O}(f,\tau )={\bf R}\backslash ({\mathcal L}(f,\tau )\bigcup 
{\mathcal N}(\tau ))\, .
$$
For each $t\in {\mathcal O}(f,\tau )$ either $|\, G(t+\frac jN\, \tau )\, |\geq 
\delta $ for all $j=0,1,\dots ,N$ or there exists a number $j_0\in \{ 0,1,\dots ,N\} 
$ such that $|\, G(t+\frac {j_0}N\, \tau )\, |< \delta $. Consider the minimal number 
$j_0$ for which the last inequality holds. If $j_0<N$, then for every $j\in \{ j_0+1,
\dots ,N\} $ we have
$$
|\, G(t+\frac jN\, \tau )-G(t+\frac {j_0}N\, \tau )\, |\geq
$$ $$
\geq |\, \Delta \sin \, \alpha (t+\frac jN\, \tau )- \Delta \sin \, \alpha (t+\frac 
{j_0}N\, \tau )\, |-|\, f(t+\frac jN\, \tau )-f(t+\frac {j_0}N\, \tau )\, |>3\delta -
\delta =2\delta \, ,
$$
and therefore, $|\, G(t+\frac jN\, \tau )\, |>\delta $. We have got that in the case 
$t\in {\mathcal O}(f,\tau )$ there is at most one number $t+\frac jN\, \tau $, $j=0,
1,\dots ,N$, such that $|\, G(t+\frac jN\, \tau )\, |<\delta $. Let $P=\{ t\in
{\bf R}:|\, G(t)\, |<\delta \} $,
$$
\widetilde {\chi}(t)\doteq \sum\limits_{j=0}^N\chi _P\, (t+\frac jN\, \tau )\, ,\ 
t\in {\bf R}\, .
$$
We have
$$
\overline {\lim\limits_{b\to +\infty}}\ \ \frac 1{2b}\ \, \int\limits_{-b}^b 
\widetilde {\chi}(t)\, dt=(N+1)\, \kappa \, (P)\, .
$$
On the other hand, $\widetilde {\chi}(t)\leq 1$ for all $t\in {\mathcal O}
(f,\tau )$, hence 
$$
\overline {\lim\limits_{b\to +\infty}}\ \ \frac 1{2b}\ \, \int\limits_{-b}^b 
\widetilde {\chi}(t)\, dt\leq
$$ $$
\leq \overline {\lim\limits_{b\to +\infty}}\ \ \frac 1{2b}\ \, \int\limits
_{[-b,b]\, \cap \, {\mathcal O}(f,\tau )}\widetilde {\chi}(t)\, dt+
\overline {\lim\limits_{b\to +\infty}}\ \ \frac 1{2b}\ \, \int\limits
_{[-b,b]\, \backslash \, {\mathcal O}(f,\tau )}\widetilde {\chi}(t)\, dt\leq
$$ $$
\leq 1+(N+1)\, \kappa \, ({\mathcal L}(f,\tau )\bigcup {\mathcal N}(\tau ))<
1+\frac 12\, (N+1)\, \epsilon \, .
$$
Therefore, $\kappa \, (P)<\frac 1{N+1}+\frac {\epsilon}2<\epsilon \, $.
\end{proof}

\begin{corollary} \label{c4.1}
Let ${\mathbb F}\in {\mathcal A}^{(B)}$, $\Delta >0$. Then for any $\epsilon \in 
(0,1]$ there exist numbers $\delta =\delta (\epsilon ,\Delta )>0$ and $\widetilde 
{\alpha}=\widetilde {\alpha}(\epsilon ,\Delta ,{\mathbb F})>0$ such that for every 
function $g\in L^{\infty}({\bf R},{\bf R})$, for which $\| g\| _{\infty}\leq \delta 
$, and for all $\alpha \geq \widetilde {\alpha}$ and all functions $f\in {\mathbb F}$
$$
\kappa \, (\{ t\in {\bf R}:|\, f(t)+\Delta \sin \alpha t+g(t)\, |<\delta \} )< 
\epsilon \, .
$$
\end{corollary}

{\it Proof of Theorem \ref{th3.1}}. Let $\Delta _0=\frac {\Delta}2\, $, $f_0(.)=
f(.)$ (for all functions $f\in {\mathbb F}$). By Corollary \ref{c4.1}, there are
numbers $\delta _0=\delta _0(\Delta )>0$ and
$\alpha _0=\alpha _0(b,\Delta ,{\mathbb F})\in \frac {2\pi }b\, {\bf N}$ such that
for all functions $f_1(t)\doteq f_0(t)+\Delta _0\sin \alpha _0t$, $t\in {\bf R}$, 
and all functions $\widetilde g_1\in L^{\infty}({\bf R},{\bf R})$, for which
$\| \widetilde g_1\| _{\infty}\leq \delta _0\, $, the inequality
$$
\kappa \, (\{ t\in {\bf R}:|\, f_1(t)+\widetilde g_1(t)\, |<\delta _0\} )<
2^{-1} \eqno (4.1)
$$
holds, furthermore $\{ f_1(.):f\in {\mathbb F}\} \in {\mathcal A}^{(B)}$. We 
shall successively for $j=1,2,\dots $ find numbers $\Delta _j=\Delta _j(\Delta )
>0$, $\delta _j=\delta _j(\Delta )>0$, 
$\alpha _j=\alpha _j(b,\Delta ,{\mathbb F})\in \frac {2\pi }b\, {\bf N}$ and
functions $f_{j+1}\in B({\bf R},{\bf R})$ dependent on $f_j\, $, $\Delta _j\, $ 
and $\alpha _j\, $, for which $\{ f_{j+1}(.):f\in {\mathbb F}\} \in {\mathcal A}
^{(B)}$. If the numbers $\Delta _k\, $, $\delta _k\, $, $\alpha _k$ and
the functions $f_{k+1}$ have been found for all $k=0,\dots ,\, j-1$, where $j\in 
{\bf N}$, then we choose the number $\Delta _j=\Delta _j(\Delta )>0$ such that 
the inequalities $\Delta _j<2^{-(j+1)}\Delta $, $\Delta _j\leq 2^{-j}\delta _0\, 
$, $\Delta _j\leq 2^{-(j-1)}\delta _1\, $, $\dots $ , $\Delta _j\leq 2^{-1}\delta 
_{j-1}\, $ hold. Further (according to Corollary \ref{c4.1}), choose numbers $\delta 
_j=\delta _j(\Delta )>0$ and $\alpha _j=\alpha 
_j(b,\Delta ,{\mathbb F})\in \frac {2\pi }b\, {\bf N}$ such that for all functions 
$f_{j+1}(t)\doteq f_j(t)+\Delta _j\sin \alpha _jt$, $t\in {\bf R}$, and all
functions $\widetilde g_{j+1}\in L^{\infty}({\bf R},{\bf R})$, for which
$\| \widetilde g_{j+1}\| _{\infty}\leq \delta _j\, $, the inequality
$$
\kappa \, (\{ t\in {\bf R}:|\, f_{j+1}(t)+\widetilde g_{j+1}(t)\, |<\delta _j\} 
<2^{-j-1} \eqno (4.2)
$$
holds. Furthermore, we also have $\{ f_{j+1}(.):f\in {\mathbb F}\} \in {\mathcal A}
^{(B)}$. Next, let us set
$$
g(t)=\sum\limits_{j=0}^{+\infty}\Delta _j\sin \, \alpha _jt\, ,\ t\in {\bf R}\, .
$$
Since $\Delta _0=\frac {\Delta}2$ and $\Delta _j<2^{-(j+1)}\Delta $ for all $j\in 
{\bf N}$, it follows that the function $g(.)$ is continuous and $b$-periodic.
Moreover,
$$
\| g\| _{\infty}\leq \sum\limits_{j=0}^{+\infty}\Delta _j<\Delta \, .
$$
We define the functions
$$
g_j(t)=\sum\limits_{k=j}^{+\infty}\Delta _k\sin \, \alpha _kt\, ,\ t\in {\bf R}\, 
,\, j\in {\bf N}\, .
$$ 
For all $t\in {\bf R}$ we have
$$
|g_j(t)|\leq \sum\limits_{k=j}^{+\infty}\Delta _k\leq \sum\limits_{k=j}^{+\infty}
2^{-k+j-1}\delta _{j-1}=\delta _{j-1}\, .
$$
Hence, it follows from (4.1) and (4.2) that for all numbers $j=0,1,\dots $ (and all
functions $f\in {\mathbb F}$) the inequality
$$
\kappa \, (\{ t\in {\bf R}:|\, f(t)+g(t)\, |<\delta _j\} <2^{-j-1}  
$$
holds. The proof of Theorem \ref{th3.1} is complete.


\begin{thebibliography}{99}

\bibitem[1]{1} A.M.~Dolbilov, I.Ya.~Shne\u iberg. {\it Almost periodic
multivalued maps and their selections}, Sibirsk. Mat. Zh., {\bf 32}:2, 172--175,
1991. English transl. in {\it Siberian Math. J.}, {\bf 32}, 1991.

\bibitem[2]{2} L.I.~Danilov. {\it Almost periodic selections of multivalued 
maps}, Izv. Otdela Mat. i Informat. Udmurtsk. Gos. Univ., Vyp.~1, Izhevsk 1993, pp.
16-78. (Russian)

\bibitem[3]{3} L.I.~Danilov. {\it On selections of multivalued almost periodic
maps}, Manuscript No.~2340--B95, deposited in VINITI 31.07.95, Moscow (1995).
(Russian)

\bibitem[4]{4} L.I.~Danilov. {\it Measure-valued almost periodic functions and
almost periodic selections of multivalued maps}, Mat. Sb., {\bf 188}:10, 3--24,
1997. English transl. in {\it Sb. Math.}, {\bf 188}, 1997.

\bibitem[5]{5} L.I.~Danilov. {\it On almost periodic multivalued maps},
Mat. Zametki, {\bf 68}:1, 82--90, 2000. English transl. in {\it Math. Notes}, 
{\bf 68}, 2000.

\bibitem[6]{6} L.I.~Danilov. {\it On equi-Weyl almost periodic selections of
multivalued maps}, Preprint arXiv: math.CA/0310010, 2003.

\bibitem[7]{7} L.I.~Danilov. {\it On Weyl almost periodic selections of
multivalued maps}, Manuscript No.~981--B2004, deposited in VINITI 09.06.2004, 
Moscow (2004). (Russian)

\bibitem[8]{8} J.~Andres. {\it Bounded, almost-periodic and periodic solutions of 
quasilinear differential inclusions}. In: "Differential Inclusions and Optimal
Control" (ed. by J.Andres, L.G\'orniewicz and P.Nistri), LN in Nonlin. Anal. {\bf
2}, 35--50, 1998.

\bibitem[9]{9} J.~Andres, A.M.~Bersani, K.~Le\'sniak. {\it On some 
almost-periodicity problems in various metrics}, Acta Appl. Math., {\bf 65}:(1-3), 
35--57, 2001.

\bibitem[10]{10} B.M.~Levitan. {\it Almost-periodic functions}, Gostekhizdat,
Moscow 1953. (Russian)

\bibitem[11]{11} J.~Marcinkiewicz. {\it Une remarque sur les espaces de 
M.~Besicowitch}, C. R. Acad. Sc. Paris. 1939. V.~208. Pp. 157-159.

\bibitem[12]{12} L.A.~Lyusternik, V.I.~Sobolev. {\it A short course in functional
analysis}, Vysshaya Shkola, Moscow 1982. (Russian)

\bibitem[13]{13} L.I.~Danilov. {\it Uniform approximation of almost periodic 
functions in the sense of Stepanov}. Izv. Vyssh. Uchebn. Zaved. Ser. Mat., 1998, 
No.~5, 10--18. English transl. in {\it Russian Math. (Iz. VUZ)}, {\bf 42}, 1998.

\bibitem[14]{14} L.I.~Danilov. {\it Uniform approximation of Stepanov almost
periodic functions and almost periodic selections of multivalued maps}, Manuscript 
No.~354--B2003, deposited in VINITI 21.02.2003, Moscow (2003). (Russian)

\bibitem[15]{15} L.I.~Danilov. {\it Uniform approximation of Stepanov almost 
periodic functions}, Izv. Inst. Mat. i Informat. Udmurtsk. Gos. Univ., Vyp.~1 (29), 
Izhevsk 2004, pp. 33-48. (Russian)

\end{thebibliography}
\end{document}